\documentclass[final,5p,times]{elsarticle}
\biboptions{numbers,sort&compress}
\usepackage{amsmath,amssymb,bbm,array,amsfonts,float,mathtools,amsthm}
\usepackage{mathrsfs}
\usepackage{graphicx}
\usepackage{subcaption}
\usepackage{booktabs}%
\usepackage{algorithm}%
\usepackage{algorithmicx}%
\usepackage[noend]{algpseudocode}%
\usepackage{listings,enumitem}%
\usepackage{xcolor}
\usepackage{todonotes}
\usepackage{comment}
\usepackage{float}
\usepackage[hyphens]{url}
\usepackage[hidelinks]{hyperref}
\hypersetup{breaklinks=true}
\usepackage{tikz}
\usetikzlibrary{arrows.meta,positioning,calc}

\newcommand{\rd}{{\rm d}}

\newcommand{\rG}{{\rm G}}

\newcommand{\R}{\mathbb{R}}
\newcommand{\C}{\mathbb{C}}

\newcommand{\SO}{{\rm SO}}

\newcommand{\SU}{{\rm SU}}

\newcommand{\U}{{\rm U}}

\renewcommand{\det}{\mathop\mathrm{det}\nolimits}

\newcommand{\End}{{\mathrm{End}}}

\renewcommand{\epsilon}{\varepsilon}

\renewcommand{\Im}{\mathop{\mathrm{Im}}}

\renewcommand{\Re}{\mathop{\mathrm{Re}}}

\newcommand{\vol}{\mathrm{vol}}

\newcommand{\qandq}{\quad\text{and}\quad}

\def\<{\mathopen{}\left<}
\def\>{\right>\mathclose{}}
\def\({\mathopen{}\left(}
\def\){\right)\mathclose{}}

\usepackage{multicol, color}

\definecolor{gold}{rgb}{0.85,.66,0}
\definecolor{cherry}{rgb}{0.9,.1,.2}
\definecolor{burgundy}{rgb}{0.8,.2,.2}
\definecolor{orangered}{rgb}{0.85,.3,0}
\definecolor{orange}{rgb}{0.85,.4,0}
\definecolor{olive}{rgb}{.45,.4,0}
\definecolor{lime}{rgb}{.6,.9,0}
\definecolor{green}{rgb}{.2,.7,0}
\definecolor{grey}{rgb}{.4,.4,.2}
\definecolor{brown}{rgb}{.4,.3,.1}

\newtheorem{theorem}{Theorem}
\newtheorem{proposition}[theorem]{Proposition}

\theoremstyle{definition}
\newtheorem{definition}[theorem]{Definition}

\begin{document}

\makeatletter
\def\ps@pprintTitle{%
  \let\@oddhead\@empty
  \let\@evenhead\@empty
  \let\@oddfoot\@empty
  \let\@evenfoot\@oddfoot}
\makeatother
\begin{frontmatter}

\title{
Neural and numerical methods for \texorpdfstring{$\mathrm{G_2}$}{G2}-structures on contact Calabi--Yau \texorpdfstring{$7$}{7}--manifolds}

\author[a]{Elli Heyes}
\author[b]{Edward Hirst}
\author[b]{Henrique N. S\'a Earp}
\author[b]{Tom\'as S. R. Silva}

\affiliation[a]{
    organization={Abdus Salam Centre for Theoretical Physics},
    addressline={Imperial College London},
    postcode={SW7 2AZ},
    country={UK}
}
\affiliation[b]{
    organization={Instituto de Matemática, Estatística e Computação Científica (IMECC) da Universidade Estadual de Campinas (UNICAMP)
    }, 
    postcode={13083-859},
    country={Brazil}  
}

\begin{abstract}
A numerical framework for approximating $\rG_2$-structure $3$-forms on contact Calabi--Yau manifolds is presented. The approach proceeds in three stages: first, existing neural network models are employed to compute an approximate Ricci-flat metric on a {Calabi–Yau} threefold. Second, using this metric and the explicit construction of a $\rG_2$-structure on the associated $7$-dimensional Calabi–Yau link in the $9$-sphere, numerical approximations of the $3$-form are generated on a large set of sampled points. Finally, a dedicated neural architecture is trained to learn the $3$-form and its induced Riemannian metric directly from data, validating the learned structure and its torsion via a numerical implementation of the exterior derivative, which may be of independent interest. 
\end{abstract}


\end{frontmatter}

\section{Introduction}
\label{sec:intro}

\subsection{Context and motivation}
Compact 7-dimensional manifolds with $\rG_2$ holonomy arise naturally in M-theory as the internal spaces needed to compactify the eleven-dimensional theory down to four dimensions in the absence of background flux. The existence of a single covariantly constant spinor on such a $\rG_2$-manifold allows for minimally supersymmetric, chiral effective theories in four dimensions. For this reason, these objects occupy a central role in attempts to connect the geometric structure of M-theory with phenomenologically realistic low-energy physics.

Compared to Calabi--Yau manifolds, which have been extensively studied and benefit from powerful tools in algebraic and complex geometry, less is known about the overall landscape of $\rG_2$-manifolds. Their classification is inherently more difficult, not least because they lack a suitable analogue for Dolbeault cohomology, which makes Calabi–Yau moduli more tractable. To date, a few systematic construction methods exist, such as Joyce's orbifold resolutions \cite{joyce1,joyce2,joyce_book}, the twisted connected sum construction, pioneered by Kovalev \cite{Kovalev} and further developed by Corti-Haskins-Nordstr\"om-Pacini \cite{HHN12,Corti:2012kd,CHNP12}, its extra-twisted generalisations such as  Nordstr\"om's extra-twisted connected sums \cite{extra_TCS}, and the Joyce-Karigiannis construction by gluing families of Eguchi–Hanson spaces \cite{joyce2021}.

These constructions often allow one to compute topological invariants such as Betti numbers, intersection data, and information about singular loci. Such topological data is useful for understanding the effective field theory arising from M-theory compactifications. Betti numbers determine the number of moduli and abelian gauge fields; singularities encode non-abelian gauge groups and chiral matter; and topological cycles control aspects of fluxes and membrane instanton effects. In this way, topology provides significant insight into the low-energy physics induced by $\rG_2$ compactifications.

However, the absence of explicit $\rG_2$ holonomy metrics remains a major obstacle. Knowing the metric would make it possible to compute geometric quantities such as cycle volumes, Yukawa couplings and the full non-perturbative potential generated by M2- and M5-brane effects. These are crucial for understanding moduli stabilisation, supersymmetry breaking, and the detailed structure of the 4d effective theory. Thus, while current constructions offer valuable topological control, obtaining or approximating the underlying metric would open the door to a much deeper understanding of M-theory compactifications on $\rG_2$-manifolds.

The past decade has seen significant progress in numerically approximating Ricci-flat Calabi–Yau metrics \cite{Donaldson3, Douglas:2020hpv, Larfors:2021pbb, Gerdes:2022nzr, Ashmore:2021ohf, Hendi:2024yin, Berglund:2022gvm, Larfors:2022nep}, initiated with computational implementation of Donaldson balanced metrics algorithm \cite{Douglas:2006hz, Douglas:2006rr}. 
One of the most important developments has been the application of machine learning and modern numerical optimisation techniques which have shown far greater efficiency than previous numerical methods.
This has enabled the first computation of physical quantities in string compactifications, such as Yukawa couplings in heterotic compactifications \cite{Butbaia:2024tje, Berglund:2024uqv, Constantin:2024yaz, Constantin:2025vyt}.
These successes provide a template for similar numerical approximations of $\rG_2$ holonomy metrics, which is the subject of this paper. 

In this work, a neural network architecture is built and trained to approximate both the defining 3-form and induced metric of the $\rG_2$-structure on a contact Calabi--Yau link $\rG_2$-manifold for the first time. 
Within this area, we note related recent work \cite{Douglas:2024pmn} which uses neural networks to approximate nowhere‐vanishing harmonic 1-forms on Calabi--Yau manifolds which feed into a Joyce‐type $\rG_2$ construction, and in later work could be used to build the holonomy metric. 
Also mentioning \cite{Brieuc2025}, which learns a Ricci-flat metric on $T^7$, which has potential to be later connected to the ``neck'' region of a twisted connected sum manifold.

\subsection{Key concepts in special geometries}
A 7-manifold $M$ is called a \emph{$\rG_2$-structure manifold} if it admits a $\rG_2$-structure, i.e. a non-degenerate 3-form $\varphi$ as defined in \S\ref{sec:g2_background}. If moreover this 3-form is both closed and coclosed with respect to the exterior derivative and the Hodge star, then $\varphi$ is said to be \emph{torsion-free}, the pair $(M,\varphi)$ is called a \emph{$\rG_2$-manifold}, and it admits a compatible metric $g_{\varphi}$ which has holonomy $\text{Hol}(g_\varphi)$ contained in the exceptional Lie group $\rG_2 \subset \SO(7)$. Finding such a holonomy $\rG_2$ metric therefore amounts to finding a torsion-free $\rG_2$-structure, which corresponds to solving a very complicated nonlinear PDE. 

In this paper, the  natural $\rG_2$-structures on certain contact Calabi--Yau (cCY) 7-manifolds is considered, which satisfy the weaker \emph{coclosed} condition $\rd\ast\varphi=0$. These manifolds are particularly interesting as they provide a natural geometric bridge between Calabi--Yau and $\rG_2$-geometries \cite{Calvo-Andrade:2016fti,Lotay_2023}, and because they are real-algebraic, hence much more amenable to computational implementation. A contact Calabi--Yau manifold carries a transverse $\SU(3)$-structure from which its canonical $\rG_2$-structure is obtained. Consequently, cCY manifolds furnish a geometrically nontrivial and numerically tractable class of $\rG_2$-structures, in contrast with the general problem of constructing torsion-free $\rG_2$-geometries where no explicit parametrisation is available. For this reason, they provide an ideal setting for numerical and machine-learning approaches to exceptional geometries.

Contact Calabi--Yau $7$-manifolds were first studied from a machine-learning perspective in \cite{Aggarwal:2023swe}, where the authors constructed an essentially exhaustive dataset of $7$-dimensional links of weighted projective Calabi--Yau $3$-fold hypersurface singularities, covering almost all possible $\mathbb{P}^4(\mathbf{w})$ weight systems. They calculated associated topological invariants -- namely the Sasakian Hodge numbers and the Crowley--Nordstr\"om invariant \cite{Biswas2010, CNInvariant}-- whose computation relies on Gr\"obner basis, and used neural networks, as well as symbolic regression, to predict these quantities directly from the ambient weight data with high accuracy. Beyond providing a substantial computational speed-up over direct Gr\"obner basis methods, the analysis also led to new conjectures concerning the topology of these manifolds and its dependence on the ambient weight systems.

\subsection{Outline}
While \cite{Aggarwal:2023swe} focused on topological invariants and their prediction via machine learning, in this paper the differential geometry of contact Calabi--Yau $7$-manifolds and the associated $\rG_2$-structures is turned to. This work is organised as follows: in \S\ref{sec:background} some background of contact Calabi--Yau manifolds and $\rG_2$-geometry is surveyed; in \S\ref{sec:data} the model for approximating the CY metric, the method of sampling data points on the CY, construction of the $\rG_2$ form, and numerical verification is described;
in \S\ref{sec:ml} the machine learning architecture of the $\rG_2$ form with results of the investigations are presented; with conclusions in \S\ref{sec:conclusion}, discussing some future prospects. Throughout this work, focus is on the Fermat quintic Calabi---Yau threefold as the base manifold for all numerical experiments. The full source code for all experiments reported in this paper is publicly available at \url{https://github.com/TomasSilva/LearningG2}\;.

\section{Background}
\label{sec:background}

\subsection{\texorpdfstring{$\mathrm{G_2}$}{G2}-geometry}
\label{sec:g2_background}

A $G$-structure on a $n$-dimensional manifold is a reduction, i.e. a principal sub-bundle, of the associated frame bundle with structure group $G$, where $G$ is some Lie subgroup of $GL(n,\mathbb{R})$. In the exceptional case where $G=\rG_2 \subseteq GL(7,\mathbb{R})$, there is a convenient way to define a $\rG_2$-structure in terms of the model exterior form  $\varphi_{0}\in\Omega^3\big(\mathbb{R}^{7}\big)$ defined by
\begin{equation}
    \varphi_{0} = e^{123} + e^{145} + e^{167} + e^{246} - e^{257} - e^{347} - e^{356},
\end{equation}
in which $(e^{1},\dots,e^{7})$ is the standard dual basis on $\big(\mathbb{R}^{7}\big)^*$, with notation $e^{ijk}\coloneqq e^{i}\wedge e^{j}\wedge e^{k}$. It is well-known that such $\varphi_0$ encodes the standard vector cross-product in $\R^7$, and it is stabilised by the action of $\rG_2$.

\begin{definition}
    A $\rG_2$-structure on $M^7$ is a smooth differential form $\varphi\in\Omega^3(M)$  such that, for every $p\in M$, there exists an oriented isomorphism $f:T_{p}M\rightarrow\mathbb{R}^{7}$ (i.e. a pointwise $\rG_2$-class of frames) satisfying $\varphi\big|_{p}=f^{*}\varphi_{0}$. 
    Using notation $Q \subset F(M)$ for the frame sub-bundle defined by the isomorphisms which send $\varphi\big|_{p}\mapsto\varphi_{0}$. It follows that $Q$ is a $\rG_2$-structure.  
\end{definition}

A $\rG_2$-structure canonically induces a Riemannian metric $g_\varphi$ and an orientation on $M^7$, and hence a Hodge star operator $*_\varphi$ on $\Omega^\bullet(M)$. 

\begin{proposition}[\cite{Karigiannis2008, karigiannis2008notesg2spin7geometry}]
\label{prop:g2_metric}
    Let $\varphi \in \Omega^{3}(M)$, and let $x^{1},\ldots,x^{7}$ be local coordinates on an open set $U \subset M$.   For $i,j \in \{1,\ldots,7\}$, define local smooth functions $B_{ij}$ by
\begin{equation}
    -6\, B_{ij}\, dx^{1} \wedge \cdots \wedge dx^{7}
    =
    \left(\frac{\partial}{\partial x^{i}} \,\lrcorner\, \varphi\right)
    \wedge
    \left(\frac{\partial}{\partial x^{j}} \,\lrcorner\, \varphi\right)
    \wedge \varphi .
\label{eq:Bij-def}
\end{equation}
    If $\varphi$ defines a $\rG_{2}$-structure with associated Riemannian metric $g_\varphi$, then
\[
B_{ij} = (g_\varphi)_{ij} \sqrt{\det g_\varphi},
\]
Consequently,
\[
\det B = (\det g_\varphi)^{\frac{9}{2}},
\qquad
\sqrt{\det g_\varphi} = (\det B)^{\frac{1}{9}},
\]
and the metric components are recovered from $\varphi$ by
\begin{equation}
\label{eq:G2_metric_dfn}
(g_\varphi)_{ij} = \frac{1}{(\det B)^{\frac{1}{9}}} \, B_{ij}.
\end{equation}
\end{proposition}

Denoting the associated Hodge-dual four-form by
\[
\psi := \ast_\varphi \varphi,
\]
the triple $(\varphi,\psi,g_\varphi)$ is such that, at every point $p\in M^7$, there exists a linear isomorphism
\[
T_p M \cong \mathbb{R}^7
\]
identifying $\varphi_p$, $\psi_p$, and $g_{\varphi,p}$ respectively with the standard $\rG_2$-invariant forms $\varphi_0$, $\psi_0$, and the Euclidean metric, written
\[
\psi_{0} = *\varphi_{0}
= e^{4567} + e^{2367} + e^{2345} + e^{1357} - e^{1346} - e^{1256} - e^{1247},
\]
\[
g_0 = (e^1)^2 + \cdots + (e^7)^2.
\]
\begin{theorem}[\cite{Fernandez1982}]
    Let $(M^7,\varphi)$ be a $\rG_2$-structure manifold. Then the following are equivalent:
    \begin{enumerate}
        \item The holonomy of $g_{\varphi}$ is contained in $\rG_2$.
        \item $\nabla\varphi=0$, where $\nabla$ is the Levi-Civita connection of $g_{\varphi}$.
        \item $\rd\varphi=0$ and $d*_{\varphi}\varphi=0$.
        \item The intrinsic torsion of the $\rG_2$-structure vanishes. 
    \end{enumerate}
\end{theorem}

Any such $g_\varphi$ is automatically Ricci-flat, and its holonomy is all of $\rG_2$ if, and only if, $\pi_{1}(M)$ is finite \cite{joyce_book}. Here, such a $\rG_2$-manifold, i.e. with finite fundamental group and admitting a torsion-free $\rG_2$-structure $\varphi$, is referred to as a \emph{$\rG_2$-holonomy manifold}. 
Therefore, one can rephrase the problem of finding a Ricci-flat $\rG_2$-holonomy metric on $M$ by finding a $\rG_2$-structure three-form $\varphi$ which is both closed ($\rd\varphi=0$) and coclosed ($\rd*_{\varphi}\varphi=\rd\psi=0$). This is a difficult system to solve since $*_{\varphi}$ depends non-linearly on $\varphi$. 
One key idea prevalent in the known methods for constructing $\rG_2$-holonomy manifolds is that one first achieves a  $\rG_2$-structure which is closed though not quite coclosed, having `small' $\rd\psi$ in some suitable norm. In this sense $\varphi$ is `approximately' torsion-free. Joyce showed \cite{joyce_book} then that there exists a nearby torsion-free $\rG_2$-structure three-form $\tilde{\varphi}=\varphi+\rd\eta$ and one can use perturbative analysis to obtain the two-form correction term $\eta$. 

This paper inaugurates a line of inquiry following a different approach to obtaining $\rG_2$- holonomy metrics, using numerical geometric flows. Considering $\rG_2$-structure three-forms on 7-dimensional contact Calabi--Yau manifolds, which are non-trivial circle fibrations over Calabi--Yau threefolds (see \S\ref{sec:cCYs} for more details). These manifolds admit $\rG_2$-structures which are are coclosed ($\rd\psi=0$) but not closed ($\rd\varphi\neq0$). One may then evolve this initial condition under e.g. the \emph{Laplacian coflow}:
\begin{equation}
   \frac{\partial \psi_{t}}{\partial t} = \Delta_{\psi_{t}}\psi_{t} := (\rd\rd^{*_{t}}+\rd^{*_{t}}\rd)\psi_{t}, 
\end{equation}
where $\psi_{t}=*_{t}\varphi_{t}$ and $*_{t}=*_{\varphi_{t}}$. 
Stationary points of this flow (if they exist) correspond to torsion-free $\rG_2$-structures.
The Laplacian coflow on contact Calabi--Yau 7-manifolds was recently studied e.g. by \cite{Lotay:2022nms,earp2024flow}, where the authors showed that the coflow reduces respectively to an ODE under a special Ansatz and to K\"ahler–Ricci-type evolution equation for the Calabi--Yau structure data on the base threefold. 

More explicitly, the following programme is proposed: (1) first construct numerically the (initial) coclosed $\rG_2$-structure three-form on a contact Calabi–Yau manifold; and then (2) use e.g. the Laplacian coflow to (potentially) evolve this structure toward a torsion-free $\rG_2$-structure three-form via Physics Informed Neural Networks (PINNs). In this paper, step (1) is addressed, with set up for step (2) in following work.

\subsection{Contact Calabi--Yau manifolds}
\label{sec:cCYs}

\begin{definition}
\label{def:cCy-def}
\label{def:contact-manifold}
    A contact manifold $(M^{2n+1},\eta,\xi)$ is a $(2n+1)$-manifold  equipped with a global smooth \emph{contact structure} $1$-form $\eta$,  which satisfies everywhere 
    \[ \eta\wedge(\rd\eta)^{n}\neq0.\] 
    
    Dual to $\eta$ one also has the unique \emph{Reeb vector field} $\xi$ defined by 
    \[ \rd\eta(\xi,\cdot)=0 \qandq \eta(\xi)=1.\] 
\end{definition}

Just as K\"ahler geometry is the intersection of symplectic and Riemannian geometry, Sasakian geometry is the natural intersection of contact and Riemannian geometry. 
One may interpret structure group reductions on an odd-dimensional contact manifold $(M^{2n+1},\eta,\xi)$ as `even-dimensional' structures `transverse' with respect to a $S^1$-action along the fibres of a submersion \[S^1\to M^{2n+1}\to V\;.\] 
In particular, Sasakian geometry may be seen as transverse K\"ahler geometry, corresponding to the reduction of the transverse holonomy group to $\U(n)$. 

\begin{definition}
    A contact Riemannian manifold $(M^{2n+1},\eta,\xi,g)$ is \emph{Sasakian} if its Riemannian cone $(M^{2n+1} \times \mathbb{R}_{>0}, t^{2}g + dt^{2})$  is K\"ahler, with K\"ahler form
    \begin{equation}
        t^{2}\,\rd\eta + 2t\,\rd t\wedge\eta.
    \end{equation}
    There thus exists a \emph{transverse} complex structure $J\in\End(TM)$, such that
    \begin{equation}
        J\circ J = -I_{TK} + \eta\otimes \xi,
    \end{equation}
    and furthermore the Reeb vector field is given by
    \begin{equation}
    \label{eq:reeb-vector-field}
        \xi = -J\left(t\frac{\partial}{\partial t}\right).
    \end{equation}
\end{definition}

One can also consider Sasakian manifolds with special transverse holonomy $\SU(n)$, in which case $M^{2n+1}$ is an $S^{1}$-bundle over a Calabi--Yau manifold $V$.

\begin{definition}
    A Sasakian manifold $(M,\eta,\xi,g,J,\Upsilon)$ is said to be a contact Calabi--Yau manifold (cCY) if $\Upsilon$ is a nowhere-vanishing transverse form of horizontal type $(n,0)$, such that 
    \begin{equation}
        \Upsilon\wedge\overline{\Upsilon} = (-1)^{\frac{n(n+2)}{2}}\omega^n
        \quad\text{and}\quad 
        \rd\Upsilon=0,
        \quad\text{with}\quad
        \omega=d\eta.
    \end{equation}
\end{definition} 

In particular when $n=3$, a cCY-structure naturally induces a coclosed $\rG_2$-structure: 
\begin{proposition}[\cite{Habib2015,Calvo-Andrade:2016fti, Gray1969}] 
\label{prop:G2_cCY}
    Every cCY manifold $(M^{7},\eta,\xi,J ,\Upsilon)$ is an $S^1$-bundle $\pi:M^{7}\to V$ over a Calabi--Yau threefold, with connection 1-form $\eta$ and curvature  $\rd\eta=\omega$,  and it carries a coclosed $\rG_2$-structure
    \begin{equation}
        \varphi :=\eta\wedge \omega +\text{Re}\Upsilon,
    \end{equation} 
    and Hodge dual $4$-form
    \begin{equation}
       \psi=\ast\varphi = \frac{1}{2}\omega\wedge\omega+ \eta\wedge\text{Im}\Upsilon, 
    \end{equation}

    \noindent with torsion $\rd\varphi= \omega\wedge\omega$ and $\rd\psi=0$.

 \noindent\textbf{N.B.:} Here differential forms and their pullbacks under $\pi$ are denoted identically.
\end{proposition}

\subsection{Calabi--Yau links}
\label{sec:links}

A polynomial $f\in\C[z_{0},\dots.,z_{n}]$ is called a weighted homogeneous polynomial of degree $d$ with weights $\textbf{w}=(w_{0},\dots,w_{n})\in\mathbb{Z}^{n+1}_{>0}$ if for any $\lambda\in\mathbb{C}^{*}$
\begin{equation}
    f(\lambda^{w_{0}}z_{0},\dots,\lambda^{w_{n}}z_{n}) = \lambda^{d}f(z_{0},\dots,z_{n}).
\end{equation}
Such an $f$ defines a non-singular affine hypersurface 
\begin{equation}
    {\mathcal{V}} := \{(z_0,\dots,z_n)\in\mathbb{C}^{n+1}\;\big|\;f(z_0,\dots,z_n)=0\}.
\end{equation}
which, in general, admits a singularity at the origin. 
\begin{definition}
\label{def:weighted-link}
    Let $f:\mathbb{C}^{n+1}\rightarrow\mathbb{C}$ be a $\textbf{w}$-weighted homogeneous polynomial with an isolated singularity at the origin. Then $K_{f}:=\mathcal{V}\cap S^{2n+1}$ is called a weighted link of degree $d$ and weight $\textbf{w}$.
\end{definition}

A weighted link $K_{f}$ of degree $d$ and weights $\textbf{w}=(w_{0},\dots,w_{n})$ is said to be a \emph{Calabi--Yau link} if
\begin{equation}
    d = \sum_{i=0}^{n} w_{i}.
\end{equation}

If $K_f$ is a $7$-dimensional weighted link, the following commutative diagram holds, where the horizontal arrows are Sasakian and K\"ahler embeddings, and the vertical arrows are principal $S^{1}$--orbibundle and orbifold Riemannian submersions, respectively.
\begin{equation}
\label{eq:commutative-diagram}
    \begin{matrix}
        K_{f} & \rightarrow & S^{9} \\
        \downarrow & & \downarrow \\
        {V} & \rightarrow & \mathbb{P}^{4}(\textbf{w})
    \end{matrix}
\end{equation}
The condition $d-\sum_{i=0}^{n}w_{i}=0$ means that the smooth weighted projective variety ${V}$ is a Calabi--Yau orbifold, thus Calabi--Yau links are non-trivial circle fibrations over Calabi--Yau $3$-orbifolds. Indeed, Milnor showed that $7$-dimensional links are $2$-connected compact smooth manifolds; and $K_f$ is the total space of a Hopf $S^1$-bundle over a (weighted) projective $3$-orbifold in $\mathbb{P}^4(\textbf{w})$~\cite{Milnor1969}.

\begin{proposition}[\cite{Habib2015}]
    Every Calabi--Yau link admits an $S^{1}$--invariant contact Calabi--Yau structure. 
\end{proposition}
\section{Data Generation \& Analysis}
\label{sec:data}

\subsection{Generating the base CY metric}
\label{sec:cymetric}
To obtain numerical samples of a coclosed $\rG_2$-structure over a Calabi--Yau link, our first step is to retrieve a numerical approximation of the K\"ahler form on a Calabi--Yau threefold serving as the base space. Throughout this work, the definiteness the base manifold was fixed to be the Fermat quintic threefold in $\C\mathbb{P}^4$, defined as the hypersurface associated to the weighted homogeneous polynomial
\begin{equation}
\label{eq:quintic-eq}
    f(z_0,\dots,z_4) = z_0^5 + \cdots + z_4^5\quad \in \C[z_0,\dots,z_4],
\end{equation}
with weights $\mathbf{w} = (1,\dots,1)$. More general weighted projective threefolds inducing Calabi-Yau links could equally well be considered.

Since the seminal work of Donaldson~\cite{Donaldson3}, a wide range of numerical techniques have been developed for approximating special Hermitian metrics on Calabi--Yau hypersurfaces; see, for example,
\cite{Ashmore:2019wzb,Douglas:2020hpv,Larfors:2021pbb,Gerdes:2022nzr,Ashmore:2021ohf,Hendi:2024yin,Berglund:2022gvm,Larfors:2022nep}. 
For the purposes of this work, the framework introduced in~\cite{Larfors:2021pbb} was adopted as the numerical starting point for this $\rG_2$ construction.

In~\cite{Larfors:2021pbb}, the authors develop a machine-learning-based library, \texttt{cymetric}, for numerically approximating Ricci-flat metrics on Calabi--Yau threefolds realized as hypersurfaces (or complete intersections). The input consists of points sampled uniformly with respect to a known measure on the CY manifold, and the output is a pointwise approximation of the Ricci-flat metric represented as a $(3,3)$ Hermitian tensor. The learning strategy is based on a family of Ans\"atze in which a neural network learns pointwise corrections to the Fubini--Study (FS) metric, denoted by $g_{\mathrm{NN}}$. The reconstructed metric is then obtained by combining the FS metric with the learned correction. In particular, in this implementation the element-wise multiplicative Ansatz
\[
    g_{\mathrm{pred}} = g_{\mathrm{FS}} + g_{\mathrm{FS}} \odot g_{\mathrm{NN}},
\]
is employed, where $\odot$ denotes the Hadamard (entrywise) product.
Furthermore, the authors of \texttt{cymetric} introduce a composite loss function designed to assess the quality of the predicted metric $g_{\mathrm{pred}}$. This loss simultaneously penalizes deviations from the Monge--Amp\`ere equation, violations of K\"ahlerity, non-vanishing Ricci curvature, and discrepancies in the total volume relative to the reference Fubini--Study metric $g_{\mathrm{FS}}$; see~\cite{Larfors:2021pbb} for further details.

\texttt{cymetric} was then trained using the Hadamard (entrywise) product Ansatz on a dataset of $200{,}000$ points sampled on the quintic, split into $90\%$ training and $10\%$ validation sets. After $300$ training epochs, the model exhibits satisfactory convergence of this composite loss function.

\subsection{Numerical checks on \texttt{cymetric} package}
\label{sec:numerics-cymetric}
Beyond the numerical convergence enforced by the training loss used in \texttt{cymetric}, an independent numerical validation of the K\"ahlerity of the predicted Hermitian metric was performed. Using the numerical exterior derivative method described in~\cite{fadel2025exteriorderivativemeanvalue}, the measure of torsion
\[
\|\rd \omega_{\mathrm{pred}}\| \;,
\]
was evaluated across the dataset.

In practice, the numerical exterior derivative algorithm $\linebreak \texttt{NED}(\alpha(p), \varepsilon)$ evaluates a differential form $\alpha$ at neighbouring points $p \pm \varepsilon e_i$, where $e_i$ are geodesic directions in the manifold chosen according to a local trivializing Euclidean chart, and reconstructs the exterior
derivative via a flux interpretation of Stokes' theorem; see~\cite{fadel2025exteriorderivativemeanvalue} for details. When applying \texttt{NED} to the predicted K\"ahler form $\omega_{\mathrm{pred}}$, whose values are provided by a neural
network rather than an analytic expression, the estimator exhibits a more sensitive
scale dependence with respect to the sampling radius $\varepsilon$.

Hence, the evaluation of $\omega_{\mathrm{pred}}(p\pm\varepsilon e_i)$ is highly sensitive to the choice of $\varepsilon$.
Consequently, as $\varepsilon$ increases the method operates in a pre-asymptotic regime in which an optimal scale must be identified in order to probe geometric information. 

\begin{figure}[t]
    \centering
    \includegraphics[width=\linewidth]{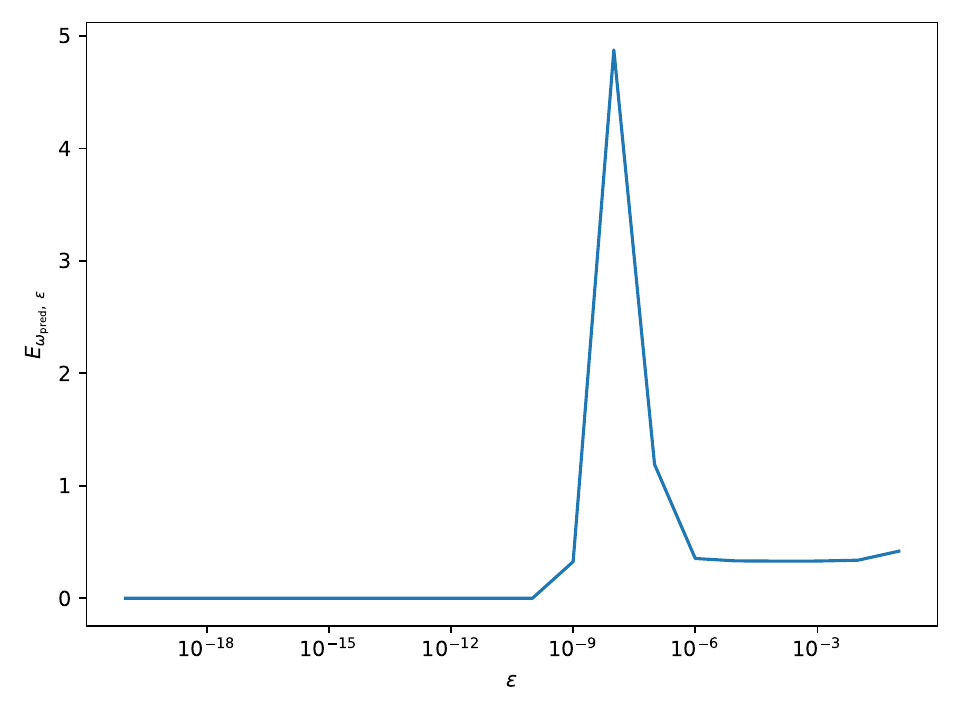}
    \caption{Variation of the K\"ahler defect with sampling radius $\varepsilon$, via \eqref{eq:medianNED}.}
    \label{fig:kahlerity_check}
\end{figure}

Figure~\ref{fig:kahlerity_check} shows the dependence of the median numerical torsion of $\omega_{\mathrm{pred}}$,
\begin{equation}\label{eq:medianNED}
    E_{\omega_{\mathrm{pred}}, \;\varepsilon} = \mathrm{median}_p
    \bigl\|
    \texttt{NED}\bigl(\omega_{\mathrm{pred}}(p),\varepsilon\bigr)
    \bigr\|\;,
\end{equation}

on the sampling radius $\varepsilon$, expressing three regimes of behaviour.
For very small $\varepsilon$, the estimator collapses to zero due to floating--point cancellation (i.e. $p+\varepsilon e_i = p-\varepsilon e_i$ up to floating point precision). 
For intermediate $\varepsilon$, a spike appears corresponding to amplification of local neural approximation errors by the discrete differential operator (i.e. where the order of $\varepsilon$ is still magnitudes below the error in the NN predictions, the division causes the \texttt{NED} to explode). 
Finally, for macroscopic $\varepsilon$ the curve stabilizes and becomes essentially independent of $\varepsilon$. 
This observed plateau is interpreted as a genuine geometric signal: at these scales the numerical exterior derivative acts as a coarse–grained operator, detecting the intrinsic torsion of the learned tensor field rather than high-frequency artefacts of the learned representation.

Ideally one would expect $\|\rd\omega_{\mathrm{pred}}\|=0$, corresponding to an exactly K\"ahler metric. 
However, since $\omega_{\mathrm{pred}}$ is produced by a neural approximation and the exterior derivative is evaluated through a discrete flux operator, the relevant diagnostic is the behaviour under refinement of $\varepsilon$. 
A discretisation artefact would either diverge as $\varepsilon\to 1$ or decay proportionally to a power of $\varepsilon$. 
Instead, a scale–independent plateau was observed, indicating the presence of a bounded intrinsic residual rather than numerical noise.
In this regime it is observed that
\[
\|\rd\omega_{\mathrm{pred}}\| = O(10^{-1}),
\]
showing that the predicted metric is approximately K\"ahler with a bounded residual error independent of the probing scale.

\subsection{Constructing \texorpdfstring{$\rG_2$}{G2}-structures numerically}
\label{sec:g2form}

Having obtained a numerical approximation of the K\"ahler form on the quintic threefold base, we proceed to the pointwise numerical construction of a coclosed $\rG_2$--structure $\varphi|_{p}$ on the associated quintic Calabi--Yau link using Proposition~\ref{prop:G2_cCY}.

From \eqref{eq:quintic-eq}, and recalling Definition~\ref{def:weighted-link}, points on the \emph{real} seven--dimensional quintic link are described by ten real coordinates $(x_0,\dots,x_4,y_0,\dots,y_4)$ satisfying
\begin{align}
    (x_0 + i y_0)^5 + \cdots + (x_4 + i y_4)^5 &= 0, \notag\\
    x_0^2 + \cdots + x_4^2 + y_0^2 + \cdots + y_4^2 &= 1.
\end{align}
In other words, $x_i=\Re(z_i)$ and $y_i=\Im(z_i)$ in the notation of \eqref{eq:quintic-eq}, and
\[
(x_0,\dots,x_4,y_0,\dots,y_4)\in \mathcal V \cap S^{9}\subset \R^{10}.
\]
Geometrically, the link is a principal $S^{1}$--bundle over the quintic Calabi--Yau threefold. 
The six real dimensions arise from the base Calabi--Yau manifold, obtained as the projective hypersurface defined by \eqref{eq:quintic-eq}, while the additional seventh dimension corresponds to the phase of the homogeneous coordinates, generating the circle fibre.

In practice, points on the link are generated from previously sampled points on the Calabi--Yau hypersurface (see \S\ref{sec:cymetric}) as follows. 
Given a point $[z_0:\cdots:z_4]$ on the projective quintic, a representative $z\in\C^5\setminus\{0\}$ satisfying the defining polynomial is first chosen, and then normalized to lie on the unit sphere $S^9\subset\C^5$. 
The $S^1$--action
\[
z \;\longmapsto\; e^{i\theta}z, \qquad \theta\in[0,2\pi),
\]
generates the circle fibre, and therefore each base point produces a one--parameter family of points on the link. 
Thus, numerically, link points are obtained by normalizing sampled Calabi--Yau points and applying a phase rotation $\theta$, providing a parametrization adapted to the $6+1$ dimensional decomposition of the geometry.

For each point $(x_0,\dots,x_4,y_0,\dots,y_4)$ on the link, the package \texttt{cymetric} is used to obtain numerical approximations of the K\"ahler form and the holomorphic volume form on the corresponding projective point $[z_0:\cdots:z_4]$ of the Calabi--Yau base.
These quantities are naturally defined on the six--(real) dimensional base manifold. They are therefore lifted to the seven--dimensional link through the projection map $\pi:K_f \to V$ appearing in the commutative diagram of \eqref{eq:commutative-diagram}. Concretely, the neural network output provides, at each sampled point, a $(3,3)$ Hermitian tensor representing the metric in complex coordinates. From this tensor the associated K\"ahler form was reconstructed as a real $(6\times 6)$ matrix in local real coordinates on the base.

The pullback $\pi^{*}\omega$ is then expressed in the link coordinate system, producing a $(7\times 7)$ real tensor whose kernel coincides with the circle fibre direction. 
An analogous construction is performed for the holomorphic volume form $\Upsilon$, which can be computed exactly on each projective representative via residues (see, e.g.,~\cite[Ch. 5]{Griffiths1994}). Its real and imaginary parts are then evaluated pointwise and subsequently pulled back to the link. Finally, the global contact $1$--form on the link is obtained by restricting to $K_f$ the dual of the Reeb vector field on $S^{9}$, given as in Definition~\ref{def:contact-manifold}.

With these components in place, the associated $\rG_2$ $3$--form $\varphi$ and its Hodge dual $4$--form $\psi$ according to Proposition~\ref{prop:G2_cCY} were constructed pointwise, represented respectively as real skew--symmetric tensors of type $(7,7,7)$ and $(7,7,7,7)$. Via Proposition~\ref{prop:g2_metric}, the $\varphi$--induced Riemannian metric $g_{\varphi}$ was also computed, represented as a real symmetric positive--definite $(7\times 7)$ matrix.
This was completed for 5 uniformly sampled $\theta$ angles per base CY point, producing a dataset of 900,000 sample points on the link; then split into 90\%:5\%:5\% for train:validation:test.

\subsection{Analysis of the numerical \texorpdfstring{$\rG_2$}{G2} dataset}
\label{sec:g2dataanalysis}
\begin{figure}[t]
    \centering
    \includegraphics[width=\linewidth]{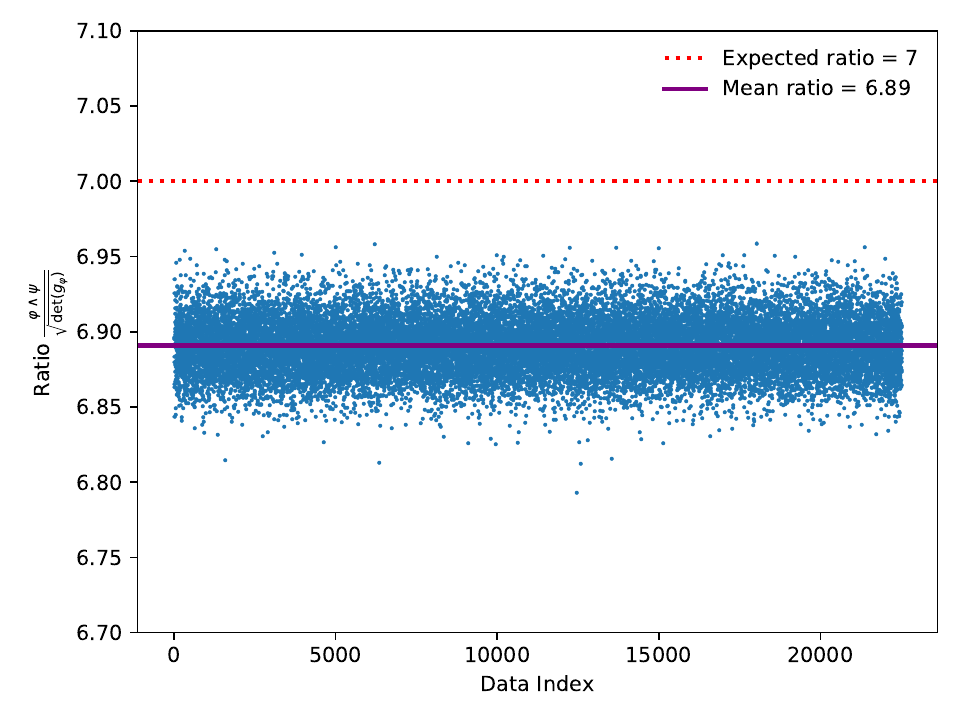}
    \caption{Pointwise values of $\frac{\varphi\wedge\psi}{\mathrm{Vol}_{g_\varphi}}$ across the dataset.}
    \label{fig:G2-identity-wedge}
\end{figure}

Beginning by assessing the correctness of the generated dataset using a fundamental $\rG_2$ algebraic identity~\cite{Karigiannis2008},
\[
\varphi \wedge \psi = 7\,\mathrm{vol}(g_\varphi),
\]
where $\mathrm{vol}(g_\varphi)$ denotes the volume form induced by the metric $g_\varphi$. In Figure~\ref{fig:G2-identity-wedge} the ratio between the computed quantity $\varphi\wedge\psi$ and $\mathrm{vol}(g_\varphi)$, is plotted pointwise over the dataset. 
The ratio is sharply concentrated around the theoretical value $7$, with only small deviations attributable to numerical discretisation and floating--point errors. This agreement provides an internal consistency check: the identity couples all components of the $\rG_2$ tensor nonlinearly and therefore tests the correctness of the entire reconstruction pipeline rather than individual tensor entries. In particular, it confirms that the induced metric $g_\varphi$ is compatible with the generated $3$--form $\varphi$ and its Hodge dual $\psi$, and that the dataset respects the defining algebraic structure of $\rG_2$ geometry. 

To gain insight into the statistical properties of the generated dataset, an exploratory data analysis of the sampled tensors was next performed. In particular, the empirical distributions of the components of $\varphi$ (Figure~\ref{fig:phi-histogram}) and the induced metric $g_{\varphi}$ (Figure~\ref{fig:g-histogram}) were examined by means of histograms computed over all sampled points. 
Figure~\ref{fig:phi-histogram} displays the empirical distributions of the $35$ oriented components of the $\rG_2$--structure $3$--form $\varphi$ across the sampled dataset. 
A striking feature is the pronounced multimodality present in several components. While some entries are sharply concentrated around zero and exhibit approximately Gaussian behaviour, many others display multiple separated peaks or asymmetric heavy tails. This behaviour is expected from the geometric origin of the data. The components of $\varphi$ are not independent scalar quantities but arise from a tensor constrained by algebraic $\rG_2$ identities and by the underlying $S^1$--fibration of the link. 

Figure~\ref{fig:g-histogram} shows the empirical distributions of the $28$ independent components of the induced metric $g_{\varphi}$ across the sampled dataset. In contrast with the $\rG_2$ three--form, the metric components exhibit a markedly different statistical behaviour. The diagonal entries are strictly positive and concentrated around well--defined nonzero values, reflecting the positive--definite character of the Riemannian metric. The off--diagonal components are centred near zero and sharply peaked, displaying rapidly decaying tails. Unlike the $\varphi$ components, the metric entries do not exhibit pronounced multimodality, indicating that $g_{\varphi}$ distributions is smoother across the manifold even with the underlying $\rG_2$ tensor transitions between different algebraic regimes.

From a learning perspective, this difference in statistics is significant: whereas the $\rG_2$ form exhibits multimodal algebraic structure, the induced metric behaves as a comparatively regular field. Consequently, predicting $g_{\varphi}$ constitutes a smoother regression problem than learning $\varphi$ itself, providing an additional diagnostic for assessing the subsequent learning performance.

\begin{figure}[t]
    \centering
    \includegraphics[width=\linewidth]{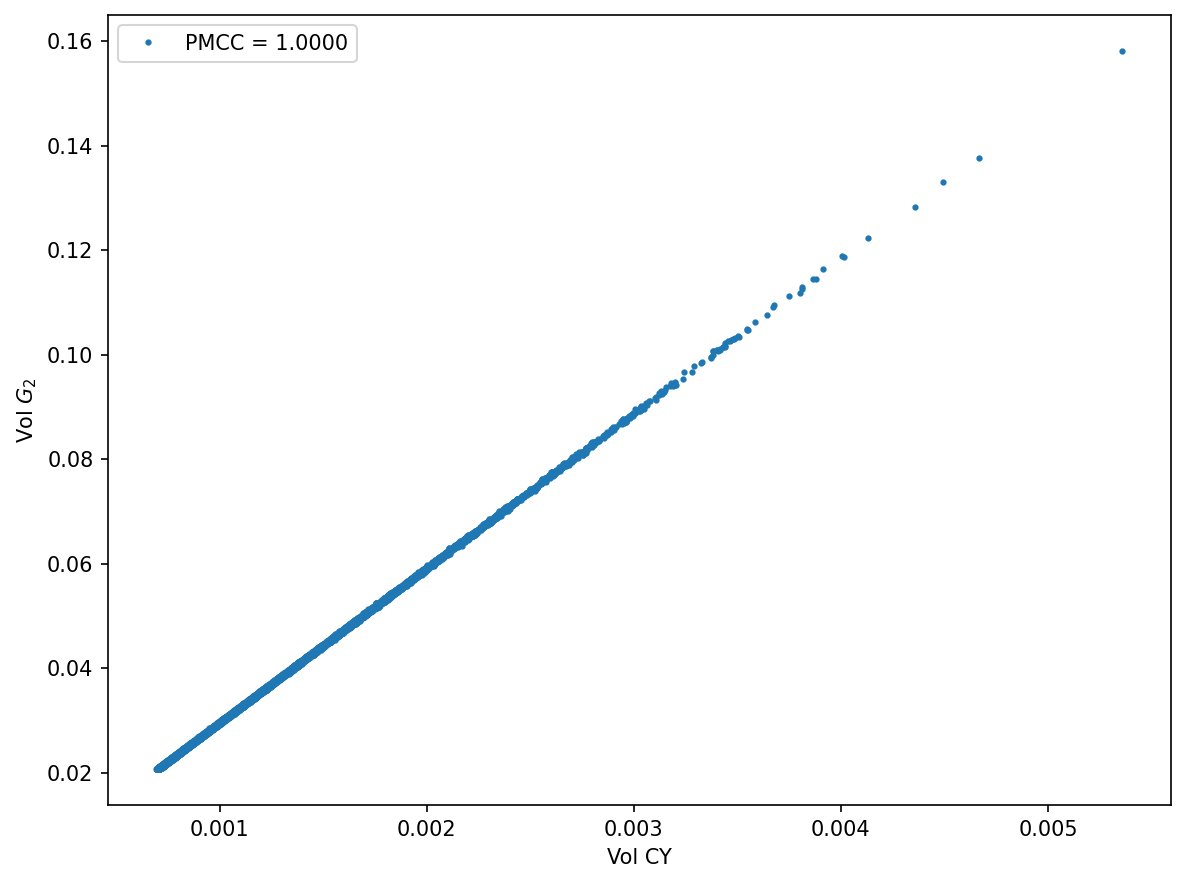}
    \caption{Comparison of the volume densities $\sqrt{\det g_{\varphi}}$ and $\sqrt{\det g_{\mathrm{CY}}}$.}
    \label{fig:vol-CY-vs-G2}
\end{figure}

Finally, in Figure~\ref{fig:vol-CY-vs-G2} the volume density induced by the $\rG_2$ metric, $\mathrm{Vol}_{g_\varphi}=\sqrt{\det g_\varphi}$, is compared with the Calabi--Yau volume density $\mathrm{Vol}_{g_{\mathrm{CY}}}$ on the base. 
The scatter plot reveals a strong linear relationship across the entire dataset, with extremely small dispersion, and a product moment correlation coefficient (PMCC) score of 1.000.
This $\mathrm{Vol}_{g_\varphi} \propto \mathrm{Vol}_{g_{\mathrm{CY}}}$ behaviour is perhaps expected from the geometric construction of the $\rG_2$ structure.
The link $7$-manifold may be considered to be composed of two parts, being either parallel or transverse to the contact form.
The transverse part correlates strongly with the CY base, and appears to be the dominant contribution to this linear relationship; whilst the parallel part parametrises the $S^1$ fibre and appears to contribute a constant factor to the $g_\varphi$ volume, shifting the linear relationship intercept above zero, as shown in Figure \ref{fig:vol-CY-vs-G2}.

This observed proportionality provides another strong global consistency check for the numerical pipeline, confirming that the reconstructed $\rG_2$ tensors satisfy the expected compatibility between the metric and the underlying complex geometry.


\section{\texorpdfstring{$\mathrm{G_2}$}{G2}-Structure Learning \& Analysis}
\label{sec:ml}

\subsection{Neural architecture}
\label{sec:architecture}
In this work, a neural architectures was designed and implemented to learn either the independent components of: 1) the $\rG_2$ defining 3-form, or 2) the $\rG_2$ metric on the quintic link, constructed for the quintic Calabi-Yau manifold with the link defining $S^1$-bundle.
The generation of this $\rG_2$ data is described in \S\ref{sec:data}. The network architectures are trained to minimise a mean squared error loss between the computed components of the $\rG_2$ data and the corresponding predictions. 

The architecture takes a composite input which has three components: \{link point ambient coordinates, $\eta$ form local coordinates, patch indices\}.
The first component is the coordinates of the considered point on the link in the ambient $\mathbb{R}^{10}$ space which the link is embedded in. The second component is the coordinates of the contact form $\eta$ in the local link coordinates in $\mathbb{R}^7$, which is the dual of the standard Reeb vector field $\xi$ on the $S^9$ restricted to the link. The third component is the two integers used to specify which patch on the base CY is used to define the local link coordinates at the input point. These components are concatenated to make the input $19$-dimensional.

The network outputs are then the independent components of either the $\rG_2$ $3$-form, or the $\rG_2$ metric in the local link coordinates. The $\rG_2$ 3-form is a $(7,7,7)$ tensor, which as a form is antisymmetric so reduces to $\binom{7}{3} = 35$ components. Conversely the $\rG_2$ metric is a $(7,7)$ tensor, which is symmetric so reduces to $7(7+1)/2 = 28$ components. This reduction to independent components both enforces the tensor symmetry structure desired and reduces the learning demand as the outputs are lower dimensional.
Therefore the NN architecture is ultimately a map of the form:
\begin{equation}
    \mathrm{NN}_\varphi: \mathbb{R}^{19} \longmapsto \mathbb{R}^{35}\;, \quad\text{or}\quad \mathrm{NN}_{g_{\varphi}}: \mathbb{R}^{19} \longmapsto \mathbb{R}^{28}\;.
\end{equation}

We remark that, for the input, using the 10-dimensional ambient coordinates for the link point is chosen to enforce a globally consistent function, as a point on the manifold has a unique description in this system, but will have multiple descriptions in local coordinates depending on the patch, which would need to be consistent (this aligns with the \texttt{cymetric} implementation also).
The $\eta$ contact form is a standard form from the Reeb vector in $S^9$ (not specific to any $\rG_2$-structure or manifold construction in general) transformed into the local coordinates on the link, it is quick to compute and empirically its inclusion in the input had a significant effect on learning results; this style of input augmentation is a standard feature engineering practise which was very useful here. 
Finally, when building the $\rG_2$ data the K\"ahler form and holomorphic volume form are provided in local coordinates on the CY base (the \texttt{cymetric} output), these two patch information indices (each a unique integer from 0 to 4) indicate which of the 5 ambient $\mathbb{C}^5 \ (\simeq \mathbb{R}^{10})$ coordinates are dropped in going to the local CY base coordinates these objects are given in; these indices also define the patching on the link which matches the CY base for simplicity. 

As demonstrated by the distribution plots given in Figures \ref{fig:phi-histogram} and \ref{fig:g-histogram}, the order of the components varies greatly. 
Therefore, to improve learning performance, within the architecture is a pair of normalisation maps.
These first fit the the training input and output data computing the mean and variance of each coordinate, then within the architecture the passed input is first passed through a normalisation layer which standardises all the coordinate orders to $O(1)$ which is optimal for learning \cite{input_norm}. 
These standardised scores are then passed to the NN architecture whose output is then passed through an equivalent inverse normalisation map on the output coordinates, such that the NN model predicts the outputs at $O(1)$ also.
Therefore the full architecture takes the inputs and predict the outputs at the correct scales, but the internal NN and the learning happen at the standardisation level, which noticeably improved learning performance.
This architecture is represented in Figure \ref{fig:architecture}. 

The internal model is a dense feed-forward NN with 4 layers of neurons with sizes (512, 512, 256, 256), with biases, and GELU activation (smoothness useful for autodifferentiation torsion computations in future work).
The Adam optimiser minimises a standard MSE loss for this supervised problem with learning rate scheduling. 
The codebase implements a Huber loss, where in practise the interpolation parameter was set to exclusively use the L2 loss (matching MSE).

\subsection{Learnt \texorpdfstring{$\rG_2$}{G2} 3-form}
\label{sec:ml_results_form}
The supervised training ran on the 900,000 link points in the train sample, validating with the 45,000 validation sample.
Training ran for 150 epochs to predict the 35 independent components of the $\varphi$ 3-form.
Loss curves through training are shown in Figure \ref{fig:loss-3form}, and the trained model was evaluated on the 45,000 test sample points producing a test MSE loss score of: $1.54 \times 10^{-6}$. 

The correlation of the model predicted 35 $\varphi$ component values, at the normalised scale, across the 45,000 test sample points against the true values are shown in Figure \ref{fig:correlation-3form}. 
Displayed with the PMCC score $\sim 1$, validating the strong learning performance.

This resulting trained model provides a smooth interpolant function between the mesh of sample points generated, and allows for quick evaluation of the $\varphi$ components at new sample points on the link.

\subsection{Learnt \texorpdfstring{$\rG_2$}{G2} metric}
\label{sec:ml_results_metric}
Equivalently, this supervised training process ran on the same 900,000 link points in the train sample, with validation on the same 45,000 validation sample points.
Training was also for 150 epochs, but now was to predict the 28 independent components of the $g_\varphi$ metric.
To ensure positive semi-definiteness of the metric this decomposition to independent components uses the Cholesky decomposition (as used in \cite{Hirst:2025seh}), via $g_\varphi = LL^T$ with $L$ a lower-triangular matrix with 28 components.
Loss curves through training are shown in Figure \ref{fig:loss-metric}, and the trained model was evaluated on the 45,000 test sample points producing a test MSE loss score of: $2.28 \times 10^{-6}$. 

The correlation of the model predicted 28 $g_\varphi$ component values, at the normalised scale, across the 45,000 test sample points against the true values are shown in Figure \ref{fig:correlation-metric}. 
Displayed with the PMCC score, $\sim 1$, again validating the strong learning performance.

This resulting trained model also provides an invaluable interpolant beyond the data sample mesh, that allows evaluation of the $\rG_2$ metric across the manifold.

\subsection{Numerical torsion verification}
\label{sec:ml_ext_deriv}
In §\ref{sec:g2form} the $\rG_2$-structure is generated pointwise on the link sample data, using the \texttt{cymetric} model to predict the $(\omega, \Upsilon)$ components and then building the $(\varphi, \psi)$ using Proposition~\ref{prop:G2_cCY}.
Since this generation is exclusively at the sampled link points, the setup is incompatible with torsion computations, where a numerical exterior derivative requires evaluating the defining $\rG_2$ forms in neighbourhoods of these points.
However, due to the results in Section \ref{sec:ml_results_form}, a continuous neural architecture model has been fitted to approximate the $\varphi$ 3-form, and can now be used to cheaply predict the values of $\varphi$ in the neighbourhoods of sampled link points as required by the numerical exterior derivative, \texttt{NED}, as in \S\ref{sec:numerics-cymetric} and~\cite{fadel2025exteriorderivativemeanvalue}.

This neural approximation to the $(\varphi, g_\varphi)$ objects now allow for computation of the $\rG_2$-structure torsion components.
To verify the $\rd\varphi$ torsion condition of this neural approximation of the $\rG_2$ structure, the trained 3-form model is used to predict the values of $\varphi$ at the 45,000 test sample link points, as well as in the neighbourhoods of each point, which the \texttt{NED} algorithm then uses to compute $\rd\varphi$ and compare to $\omega \wedge \omega$ at each sample point, using the predicted values of $\omega$ coming from the trained \texttt{cymetric} model.

Equivalently, the $\rd\psi$ torsion condition can be verified by using the trained 3-form model to again predict the values of $\varphi$ at the test sample points and their neighbourhoods, then also using the trained metric model to predict $g_\varphi$ at these same points.
This allows computation of $\psi$ at these test points and their neighbourhoods using the $\psi = \ast_{\varphi}\varphi$ relation; where this Hodge star operation requires the $g_\varphi$ metric.

The results of these torsion checks are represented across the test sample of the link manifold in Figure \ref{fig:torsion_checks} (taking $\varepsilon=10^{-5}$ in \texttt{NED} algorithm).
As well as plots of the torsion condition over the manifold test samples, a plot of the $\varphi \wedge \psi = 7\,\mathrm{vol}(g_\varphi)$ condition is also given, as the equivalent of Figure~\ref{fig:G2-identity-wedge} with all $(\varphi, \psi, g_\varphi)$ components now built from the trained NN models on these unseen test link samples, and validates expected behaviour with a ratio $\sim 7$.
For the torsion conditions, the mean average values of these checks, with standard deviation, were
\begin{equation}
    \frac{||\rd\varphi||}{||\omega \wedge \omega||} = 2.77 \pm 0.98\;,\quad ||\rd\psi|| = 0.074 \pm 0.036\;,
\end{equation} 
and additionally the full MSE between $\rd\varphi$ and $\omega \wedge \omega$ was: 2.22. 

For different choices of $\varepsilon$ in the \texttt{NED} algorithm, the qualitative behaviour is similar to that discussed in \S\ref{sec:numerics-cymetric}. 
Ideally one would expect $\rd\psi = 0$ and $\rd\varphi = \omega\wedge\omega$; however, due to the accumulation of discretisation and neural approximation errors throughout the pipeline, the observed \texttt{NED}'s behaviour instead exhibits a scale-dependent regime followed by a stable residual plateau. 
In this regime it is observed that
\[
\|\rd\psi\| = O(10^{-2}) \qquad\text{and}\qquad \|\rd\varphi\| = O(\|\omega\wedge\omega\|).
\]

These $\varepsilon$-asymptotic boundedness properties are consistent with the torsion behaviour expected from Proposition~\ref{prop:G2_cCY}, and confirm that these neural architectures successfully approximate the $\rG_2$-structure on a CY link for the first time.


\section{Conclusion}
\label{sec:conclusion}

Manifolds endowed with $\rG_2$-structures occupy a central position at the interface of differential geometry and M--theory, yet explicit analytic control over such geometries remains extremely limited. The development of numerical and data–driven approaches capable of probing these geometric structures therefore provide a complementary avenue for investigating their geometric properties and applications beyond the reach of current analytic techniques.

In this work a numerical dataset of coclosed $\rG_2$-structures arising from contact Calabi--Yau $7$-manifolds was constructed and a neural network regressor was trained to learn the associated geometric tensors, achieving mean squared errors of $<0.001$. The consistency of this dataset, as well as the trained NN models to approximate the $(\varphi, g_\varphi)$ objects, was validated by verifying $\rG_2$ compatibility identities, and showed that the torsion components of the predicted structures were numerically consistent with the theoretical expectations.

Beyond the immediate approximation problem, the present framework opens several directions for future investigation. First, the construction may be extended to other Calabi--Yau bases, allowing a systematic exploration of broader classes of $\rG_2$-structures. Second, the learned metrics provide natural initial data for neural metric flows  \cite{Halverson:2023ndu}, where one may attempt to minimise the non–vanishing torsion component and evolve towards torsion–free configurations. Finally, the approach offers a concrete setting for the numerical study of geometric evolution equations associated to $\rG_2$-structures, including Laplacian and related flows~\cite{Lotay:2022nms,earp2024flow}. We expect that combining geometric structure with machine learning optimisation will provide a practical route towards constructing and analysing novel metrics of exceptional holonomy.

\newpage
\section*{Acknowledgements}
The authors would like to thank Prof. Yang-Hui He for helpful comments in the early stages of this project; as well as Prof. Magdalena Larfors and Prof. Fabian Ruehle for discussions about optimal use of their \texttt{cymetric} package \cite{Larfors:2021pbb}.

\noindent\textbf{E.~Heyes} is supported by Schmidt Sciences AI in Science fellowship. \\
\textbf{E.~Hirst} is supported by São Paulo Research Foundation (FAPESP) grant 2024/18994-7. \\
\textbf{H.~Sá Earp} is supported by grants:\\
\makebox[5pt]{} -- FAPESP 2020/09838-0,\\
\makebox[12pt]{} with BI0S: Brazilian Institute of Data Science.\\
\makebox[5pt]{} -- FAPESP 2021/04065-6,\\
\makebox[12pt]{} with BRIDGES: Brazil-France interplays in Gauge Theory, \\ \makebox[12pt]{} extremal structures and stability.\\
\makebox[5pt]{} -- FAPESP 2024/00923-6,\\
\makebox[12pt]{} with CBG: Brazilian Centre for Geometry.\\
\makebox[5pt]{} -- 307145/2025-5 level PQ-A,\\
\makebox[12pt]{} with the Brazilian National Council for Scientific and \\ \makebox[12pt]{} Technological Development (CNPq).\\
\textbf{T.~Silva} is supported by FAPESP grant 2022/09891-4.

\noindent This research utilised computational resources of the ``\textit{Centro Nacional de Processamento de Alto Desempenho em São Paulo} (CENAPAD-SP)".

\bibliographystyle{elsarticle-num} 
\bibliography{references.bib}

\appendix
\section{Data \& Learning Plots}
\label{sec:plots}
\noindent Below are supporting plots and diagrams that demonstrate the multimodal nature of the datasets, diagrammatic representation of the NN architecture, training loss curves, model prediction correlations, and the predicted $\rG_2$ torsion components.

\begin{figure*}[!ht]
    \centering
    \includegraphics[width=\textwidth]{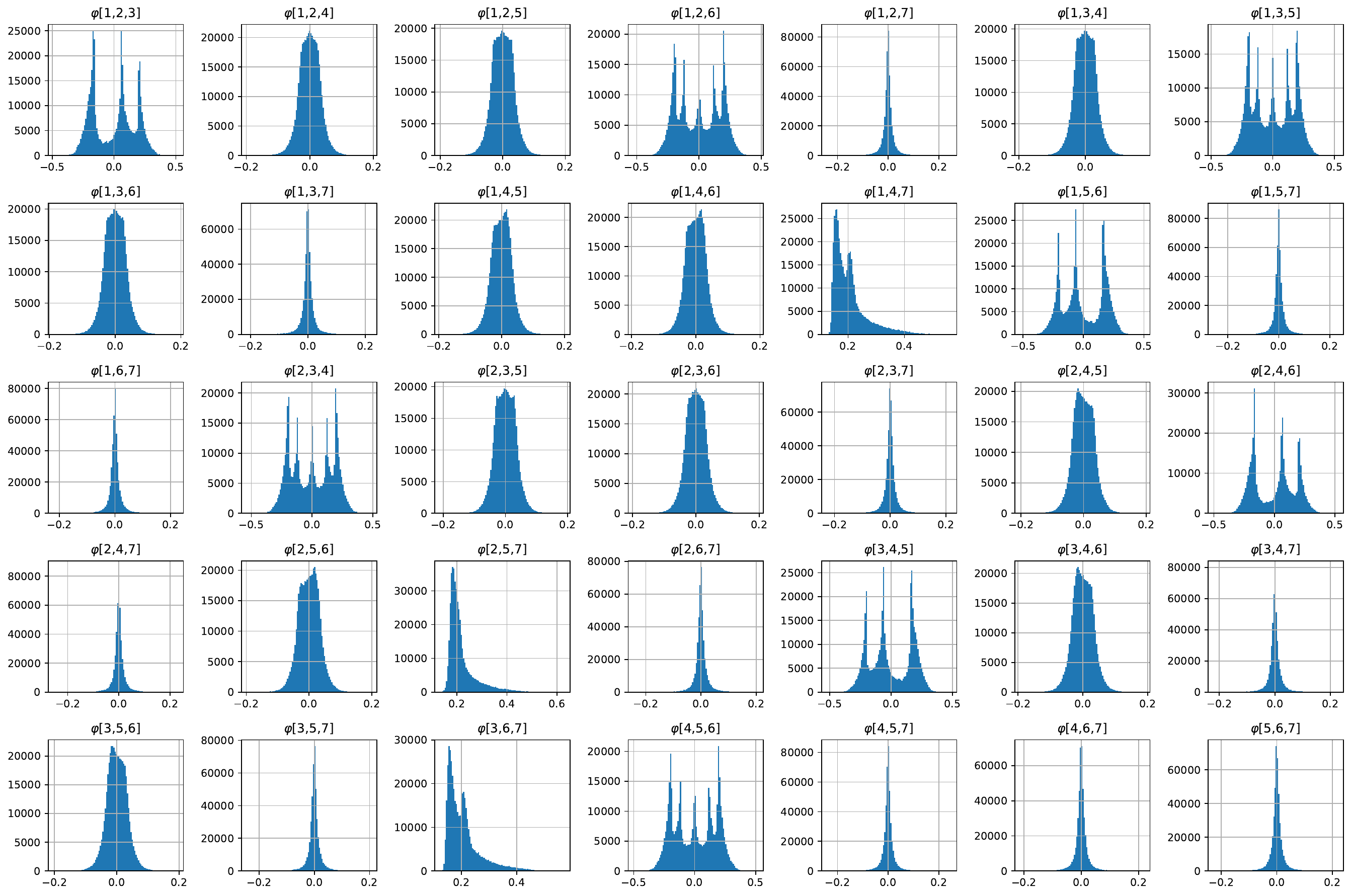}
    \caption{Histogram of the 35 independent components of the $\mathrm{G}_2$ three-form $\varphi$.}
    \label{fig:phi-histogram}
\end{figure*}

\begin{figure*}[!ht]
    \centering
    \includegraphics[width=\textwidth]{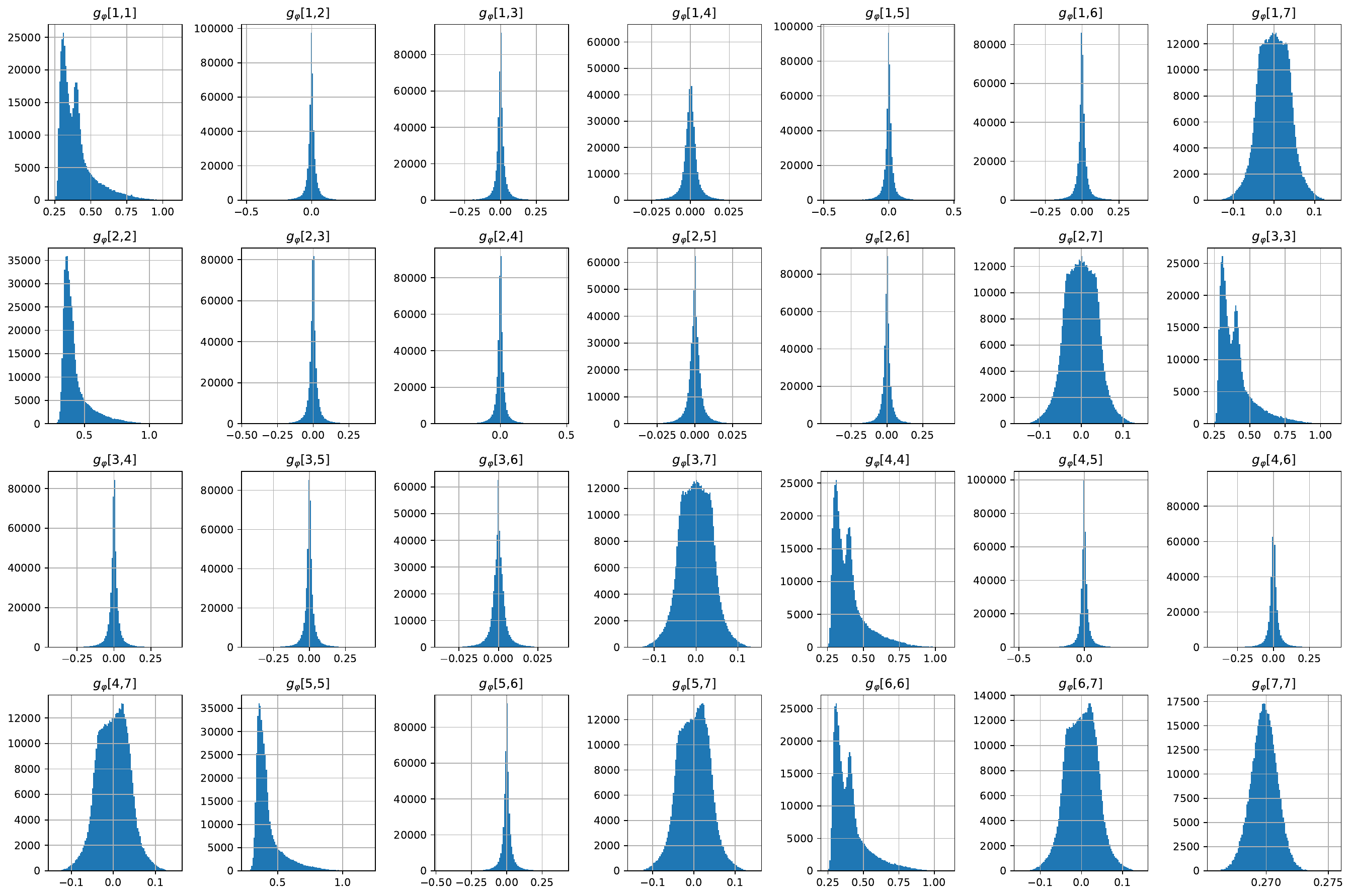}
    \caption{Histogram of the 28 independent components of the $\mathrm{G}_2$ metric $g_{\varphi}$.}
    \label{fig:g-histogram}
\end{figure*}

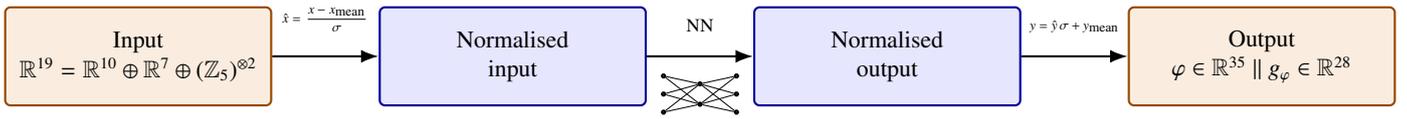
\begin{figure*}[!ht]
\centering

\begin{tikzpicture}[
  font=\small,
  >=Latex,
  box/.style={draw, rounded corners=2pt, thick, align=center, minimum width=35mm, minimum height=13mm},
  outer/.style={box, fill=orange!12, draw=orange!70!black},
  inner/.style={box, fill=blue!10, draw=blue!60!black},
  arrow/.style={-Latex, thick}
]

\node[outer] (b1) {Input\\$\mathbb{R}^{19} = \mathbb{R}^{10} \oplus \mathbb{R}^{7} \oplus (\mathbb{Z}_{5})^{\otimes 2}$};
\node[inner, right=14mm of b1] (b2) {Normalised\\input};
\node[inner, right=14mm of b2] (b3) {Normalised\\output};
\node[outer, right=14mm of b3] (b4) {Output\\$\varphi \in \mathbb{R}^{35}\ \|\ g_\varphi \in \mathbb{R}^{28}$};

\draw[arrow] (b1) -- (b2);
\draw[arrow] (b2) -- (b3);
\draw[arrow] (b3) -- (b4);

\node[above=2mm of $(b1)!0.5!(b2)$, font=\tiny] {$\hat{x}=\dfrac{x-x_{\text{mean}}}{\sigma}$};
\node[above=2mm of $(b2)!0.5!(b3)$, font=\scriptsize] {NN};
\node[above=2mm of $(b3)!0.5!(b4)$, font=\tiny] {$y=\hat{y}\,\sigma+y_{\text{mean}}$};

\begin{scope}[shift={($(b2)!0.5!(b3)+(0,-5mm)$)}, scale=0.30]
  \tikzset{nnode/.style={circle, draw, thick, fill=white, inner sep=0.40pt}}
  \foreach \i/\y in {1/-0.8,2/0,3/0.8} {\node[nnode] (L\i) at (-1.6,\y) {};}
  \foreach \j/\y in {1/-0.45,2/0.45} {\node[nnode] (M\j) at (0,\y) {};}
  \foreach \k/\y in {1/-0.8,2/0,3/0.8} {\node[nnode] (R\k) at (1.6,\y) {};}
  \foreach \i in {1,2,3} {
    \foreach \j in {1,2} {\draw[line width=0.22pt] (L\i) -- (M\j);} }
  \foreach \j in {1,2} {
    \foreach \k in {1,2,3} {\draw[line width=0.22pt] (M\j) -- (R\k);} }
\end{scope}

\end{tikzpicture}
\caption{Diagrammatic representation of the implemented architecture.}
\label{fig:architecture}
\end{figure*}

\begin{figure*}[!ht]
    \centering
    \begin{subfigure}[t]{0.48\textwidth}
        \centering
        \includegraphics[width=0.9\textwidth]{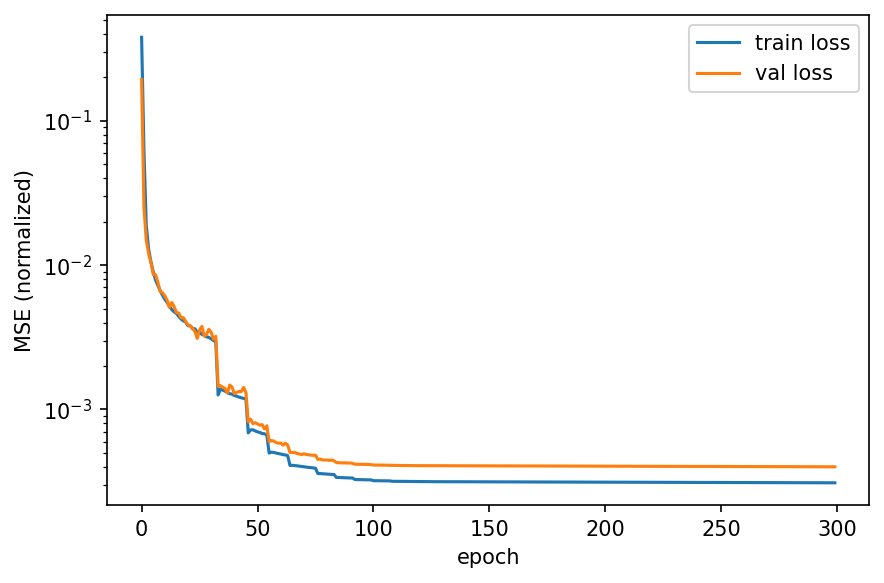}
        \caption{}
        \label{fig:loss-3form}
    \end{subfigure}
    \hfill
    \begin{subfigure}[t]{0.48\textwidth}
        \centering
        \includegraphics[width=0.9\linewidth]{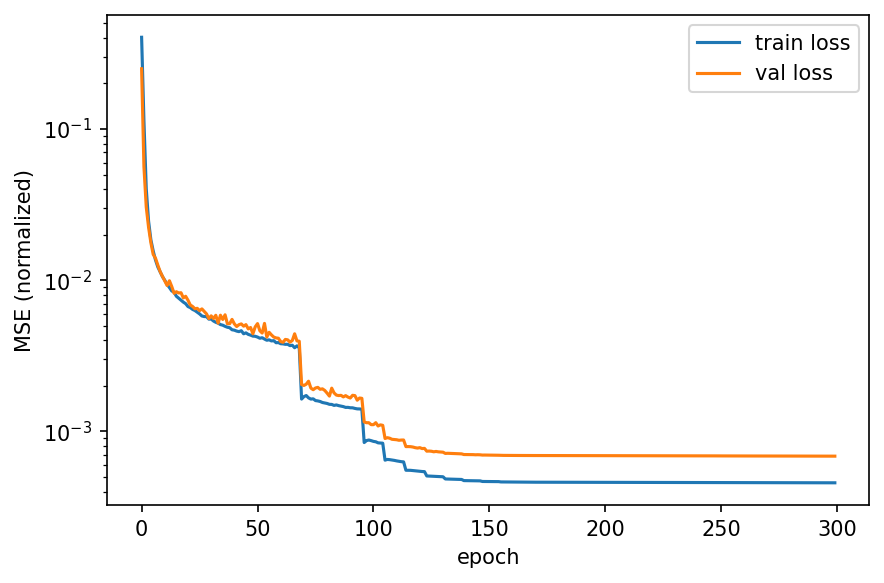}
        \caption{}
        \label{fig:loss-metric}
    \end{subfigure}
    \caption{Loss curves for the supervised training of the (a) 3-form $\varphi$ model, and (b) the $\rG_2$ metric $g_\varphi$ model.}
    \label{fig:loss_curves}
\end{figure*}

\begin{figure*}[!ht]
    \centering
    \begin{subfigure}[t]{0.48\textwidth}
        \centering
        \includegraphics[width=0.9\textwidth]{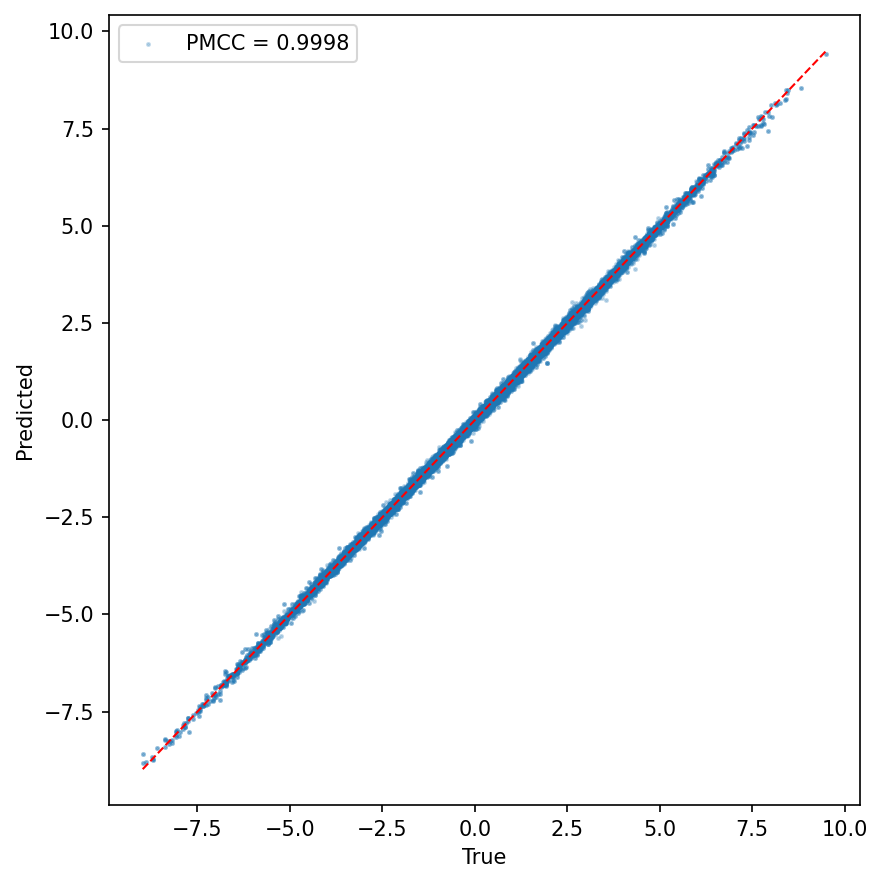}
        \caption{}
        \label{fig:correlation-3form}
    \end{subfigure}
    \hfill
    \begin{subfigure}[t]{0.48\textwidth}
        \centering
        \includegraphics[width=0.9\textwidth]{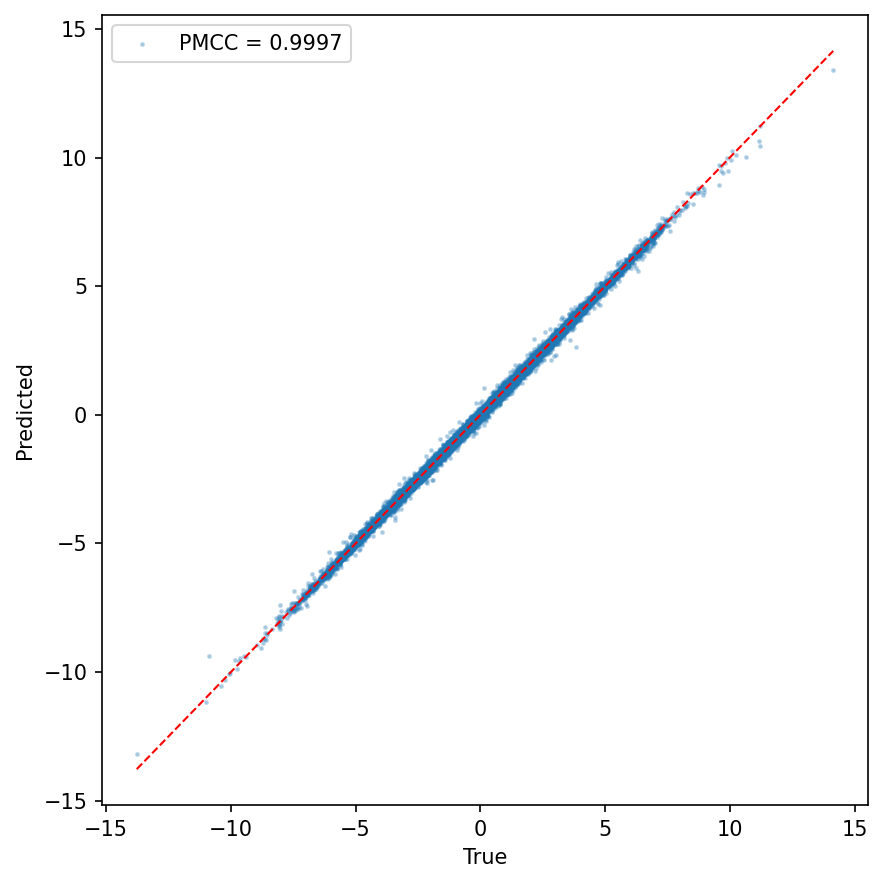}
        \caption{}
        \label{fig:correlation-metric}
    \end{subfigure}
    \caption{Correlations of the true output components against the model predicted values for the (a) 3-form $\varphi$, and (b) the $\rG_2$ metric $g_\varphi$ models; with PMCC scores.}
    \label{fig:correlations}
\end{figure*}

\begin{figure*}[!ht]
    \centering
    \begin{subfigure}[t]{0.48\textwidth}
        \centering
        \includegraphics[width=0.95\textwidth]{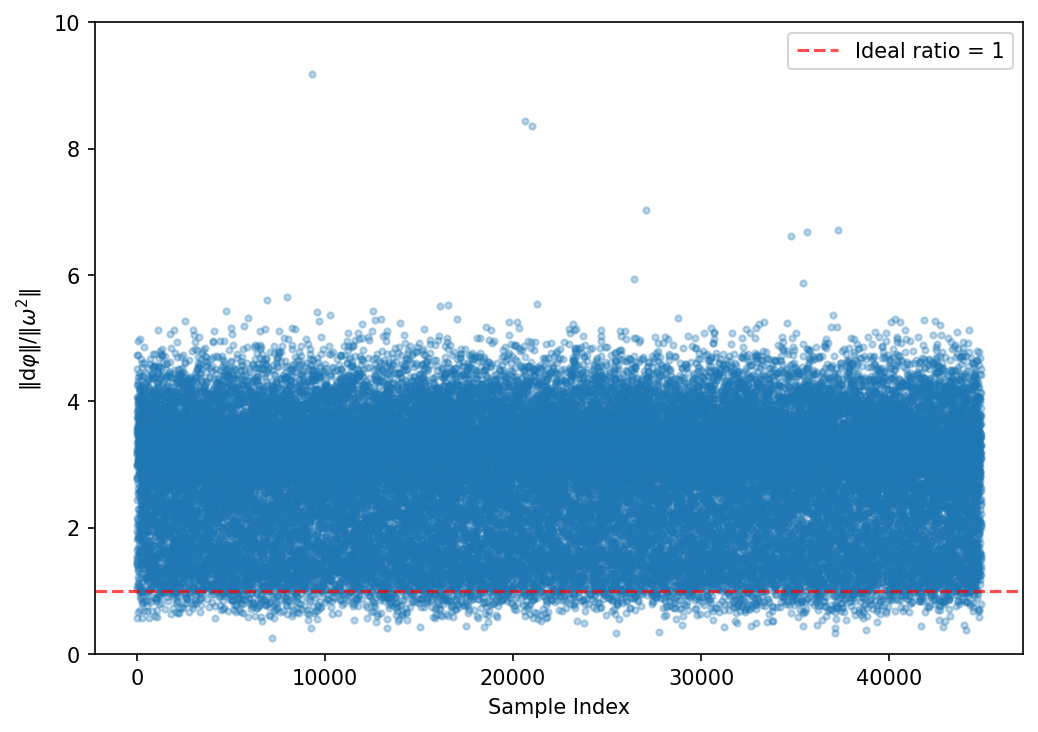}
        \caption{}
        \label{fig:torsion-phi}
    \end{subfigure}
    \hfill
    \begin{subfigure}[t]{0.48\textwidth}
        \centering
        \includegraphics[width=0.95\textwidth]{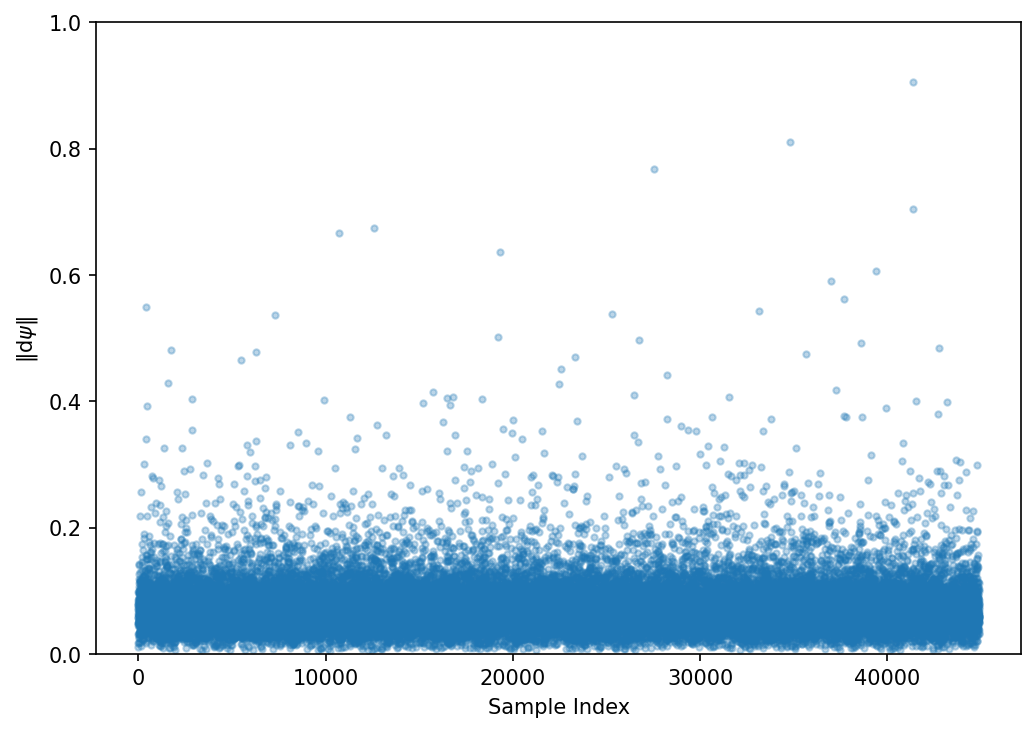}
        \caption{}
        \label{fig:torsion-psi}
    \end{subfigure}
    \\
    \begin{subfigure}[t]{0.6\textwidth}
        \centering
        \includegraphics[width=0.9\textwidth]{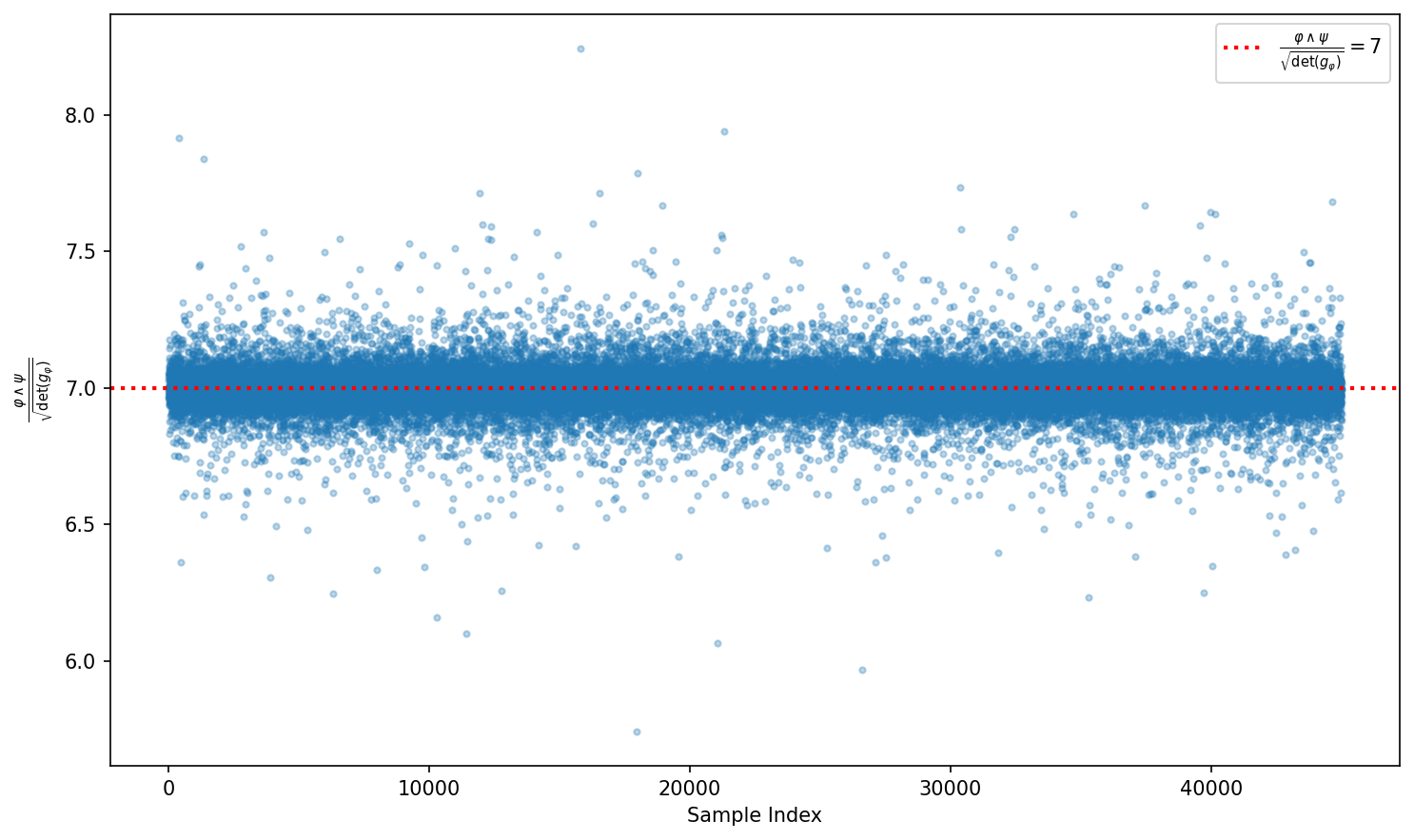}
        \caption{}
        \label{fig:torsion-vol}
    \end{subfigure}
    \caption{Torsion components of the NN learnt $\rG_2$-structure. The $\rd\varphi = \omega \wedge \omega$ torsion condition is represented by points in (a) $\sim 1$, whilst the $\rd\psi=0$ torsion condition is represented by points in (b) $\sim 0$. Additionally, the $\varphi \wedge \psi = 7 \vol (g_\varphi)$ condition for all ($\varphi, g_\varphi, \psi$) components predicted, is represented by points in (c) $\sim 7$.}
    \label{fig:torsion_checks}
\end{figure*}

\end{document}